\documentclass[reqno]{amsart}
\usepackage{amssymb,amsmath,amsthm,epsfig}
\usepackage{hyperref}

\theoremstyle{plain}
\newtheorem{thm}{Theorem}
\newtheorem{prop}[thm]{Proposition}

\newtheorem{conj}[thm]{Conjecture}
\newtheorem{cor}[thm]{Corollary}

\theoremstyle{definition}

\newcommand{\0}{\emptyset}
\newcommand{\Case}[2]{\smallskip\noindent\textbf{Case~{#1}:}~{#2}.}
\newcommand{\commentout}[1]{}
\newcommand{\defterm}[1]{\emph{#1}}
\newcommand{\diam}{\mathop{\rm diam}\nolimits}
\newcommand{\ep}{\varepsilon}
\newcommand{\isom}{\equiv}
\newcommand{\partn}{\vdash}
\newcommand{\scp}[2]{\big\langle{#1},\;{#2}\big\rangle}

\newcommand{\st}{~|~}
\newcommand{\type}{\mathop{\rm type}\nolimits}

\newcommand{\CSF}{{\mathbf X}}
\newcommand{\conn}{{\mathbf K}}
\newcommand{\subt}{{\mathbf S}}
\newcommand{\TPCSF}{{\tilde{{\mathbf X}}}} % two-part portion

\newcommand{\Pp}{{\mathbb P}}
\newcommand{\Qq}{{\mathbb Q}}

\begin{document}

\title[Chromatic symmetric functions]{On distinguishing trees by their chromatic symmetric functions}

\author{Jeremy L.\ Martin, Matthew Morin, and Jennifer D.\ Wagner}

\address{%Jeremy L.\ Martin,
Department of Mathematics,
University of Kansas,
405~Snow Hall,
1460~Jayhawk Blvd.,
Lawrence,~KS 66045, USA}
\email{jmartin@math.ku.edu}

\address{%Matthew Morin,
Department of Mathematics,
University of British~Columbia,
Room~121, 1984~Mathematics~Road,
Vancouver,~BC, Canada, V6T~1Z2}
\email{mjmorin@math.ubc.ca}

\address{%Jennifer D.\ Wagner,
Department of~Mathematics and~Statistics,
Washburn~University,
1700~SW~College~Ave.,
Topeka,~KS 66621, USA}
\email{jennifer.wagner1@washburn.edu}

\thanks{The first author was supported in part by the University of Kansas
New Faculty General Research Fund.  The second author
was supported in part by the National Sciences and Engineering Research Council of Canada.}

\begin{abstract}

Let $T$ be an unrooted tree. The \emph{chromatic symmetric function} $\CSF_T$, introduced
by Stanley, is a 
sum of monomial symmetric functions corresponding to proper colorings of $T$.  The 
\emph{subtree polynomial} $\subt_T$, first considered under a different name by
Chaudhary and Gordon, is the bivariate generating function for subtrees of $T$ by 
their numbers of edges and leaves. We prove that
$\subt_T = \langle\Phi,\CSF_T\rangle$, where $\langle\cdot,\cdot\rangle$ is the
Hall inner product on symmetric functions and $\Phi$ is a certain symmetric function
that does not depend on $T$.  Thus the 
chromatic symmetric function is a stronger isomorphism invariant than the subtree polynomial.
As a corollary, the path and degree sequences of a tree can be obtained from its chromatic 
symmetric function.  As another application, we exhibit two infinite families of trees
(\emph{spiders} and some \emph{caterpillars}), and one family of unicyclic graphs (\emph{squids})
whose members are determined completely by their chromatic symmetric functions.

\end{abstract}

\maketitle

%=================================================================
\section*{Introduction}

Let $G$ be a simple graph with vertices $V(G)$ and edges $E(G)$, and let
$\Pp$ denote the positive integers.  A \emph{(proper) coloring}
of $G$ is a function $\kappa:V(G)\to\Pp$ such that $\kappa(v)\neq\kappa(w)$ whenever
the vertices $v,w$ are adjacent.  Stanley (\cite{Stanley}; see also
\cite[pp.~462--464]{ECII}) defined the \emph{chromatic 
symmetric function} of $G$ as
  $$\CSF_G = \CSF_G(x_1,x_2,\dots) = \sum_\kappa \prod_{v\in V(G)} x_{\kappa(v)},$$
the sum over all colorings $\kappa$, where
$x_1,x_2,\dots$ are  countably infinitely many commuting indeterminates.
This definition is invariant under permutations of $\{x_i\}$, so
$\CSF_G$ is a symmetric function, homogeneous of degree $n=\#V(G)$.

The chromatic symmetric function is a much stronger isomorphism invariant than the
well-known \emph{chromatic polynomial} $\chi_G(k)$, a polynomial function of~$k$
that gives the number of colorings of
$G$ using at most $k$ colors.  Indeed, for any integer~$k$, the number
$\chi_G(k)$ can be obtained from $\CSF_G$
by setting $x_1=\cdots=x_k=1$ and $x_i=0$ for all $i>k$.

\begin{figure}
\begin{center}
\resizebox{2.0in}{0.4in}{\includegraphics{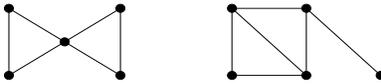}}
\end{center}
\caption{\label{sameXG-figure}
         Stanley's example \cite{Stanley} of two non-isomorphic graphs with the
         same chromatic symmetric function.}
\end{figure}

It is natural to ask whether $\CSF_G$ is a complete isomorphism invariant; that is,
whether two non-isomorphic graphs must have different chromatic symmetric functions.
The answer is no; the smallest example, shown in Figure~\ref{sameXG-figure}, was given
by Stanley in \cite{Stanley}.  Brylawski \cite{Brylawski} introduced a graph invariant
called the \emph{polychromate} and constructed an infinite family of pairs
of nonisomorphic graphs with the same polychromate.  Sarmiento \cite{Sarmiento}
proved that the polychromate is equivalent to the \emph{$U$-polynomial} studied
by Noble and Welsh \cite{NobleWelsh}, a stronger invariant than $X_G$;
therefore, each pair in Brylawski's construction
shares the same chromatic symmetric functions.  (We thank Anna de~Mier
for directing our attention to these results.)  In another direction, the non-commutative
version of the chromatic symmetric function, studied by Gebhard and Sagan \cite{GS}, is
easily seen to be a complete invariant.

Stanley's question of whether $\CSF_G$ is a complete isomorphism 
invariant for \emph{trees} remains open.  This is equivalent to the problem
of whether a tree is determined by its $U$-polynomial, since the formula
for $\CSF_G$ in terms of the $U$-polynomial \cite[Theorem~6.1]{NobleWelsh}
is easily seen to be reversible for trees.  Stanley's question was answered
in the affirmative for certain special kinds of trees
by Fougere~\cite{Fougere} and the second author~\cite{Morin}, both of
whom listed several other tree invariants that can be extracted from
the chromatic symmetric function.  Additionally,
Tan~\cite{Tan} has verified computationally
that the answer is ``yes'' for trees with 23 or fewer vertices.  (In contrast, the
chromatic polynomial is nearly useless for distinguishing trees,
because $\chi_T(k)=k(k-1)^{n-1}$ for every tree~$T$ with~$n$ vertices.)

Our main tool is Stanley's expansion of the chromatic symmetric function
in the basis of power-sum symmetric functions $p_\lambda$ \cite[Theorem~2.5]{Stanley};
see equation \eqref{Stanley's-expansion} below.
When $T$ is a (possibly trivial) tree, the coefficient $c_\lambda(T)$ of $p_\lambda$ in $\CSF_T$
has a particularly simple combinatorial interpretation.  For $A\subseteq E(T)$,
define the \emph{type} of $A$ to be the partition whose parts are the
sizes of the vertex sets of the graph with vertices $V(T)$ and edges $A$
(see Figure~\ref{type-figure} for an example).  Then,
up to a sign, $c_\lambda(T)$ is the number of edge sets $A$ of
type~$\lambda$.  As we will see, many other invariants of~$T$
can be recovered from $\CSF_T$.

Recall that the \emph{degree} of a vertex is the number of edges incident to it.
A \emph{leaf} of a tree is a vertex of degree~1, and the
unique incident edge is called a \emph{leaf edge}.
We define the \emph{subtree polynomial} of $T$ by
  $$\subt_T = \subt_T(q,r) = \sum_{\text{subtrees } S} q^{\#S} r^{\#L(S)},$$
where the sum runs over all subtrees $S$ of~$T$ with at least one edge,
and $L(S)$ denotes the set of leaf edges of~$S$.  Setting $q=t(z+1)$
and $r=1/(z+1)$ in $\subt_T$ recovers the polynomial $f_E(T;t,z)$ studied
by Chaudhary and Gordon in \cite[Section~3]{CG}.  Conversely,
$f_E(T;qr,(1-r)/r)=\subt_T(q,r)$, so the two polynomials provide identical
information about $T$.

For every non-empty set $A\subseteq E(T)$, there is a unique minimal subset
$K(A)\subseteq E(T)-A$,
called the \emph{connector} of $A$, such that $A\cup K(A)$ is a tree.
(So $K(A)=\0$ if and only if $A$ is itself a subtree of~$T$.)
The \emph{connector polynomial} of $T$ is then defined as
  $$\conn_T = \conn_T(x,y) = \sum_{\0\neq A\subseteq E(T)} x^{\#A} y^{\#K(A)}.$$

The polynomials $\subt_T$ and $\conn_T$ provide equivalent information about $T$;
we will prove in Proposition~\ref{subt-conn-equivalence} below that each of these invariants can
be obtained from the other.  Moreover, the path sequence and
degree sequence of $T$ can easily be recovered from $\subt_T(q,r)$,
as observed by Chaudhary and Gordon
\cite[Proposition~18]{CG}.

For a partition $\lambda=(\lambda_1\geq\lambda_2\geq\cdots\geq\lambda_\ell)\partn n$,
and integers $a,b,i,j$, define
  \begin{align}
  \psi(\lambda,a,b) &= (-1)^{a+b} \binom{\ell-1}{\ell-n+a+b} \sum_{k=1}^{\ell} \binom{\lambda_k-1}{a},
    \label{define-psi}\\
  \phi(\lambda,i,j) &= (-1)^{i+j} \binom{\ell-1}{\ell-n+i} \sum_{d=1}^j (-1)^d \binom{i-d}{j-d}
    \sum_{k=1}^\ell \binom{\lambda_k-1}{d}.
    \label{define-phi}
  \end{align}
We can now state our main theorem, which asserts that the subtree and connector polynomials
can be recovered from the chromatic symmetric function $\CSF_T$.

\begin{thm} \label{main-thm}
For every $n\geq 1$, and for every tree $T$ with $n$ vertices,
  \begin{equation} \label{conn-from-XT}
  \conn_T(x,y) ~=~ \sum_{a>0} \sum_{b\geq 0} x^a y^b
  \sum_{\lambda\partn n} \psi(\lambda,a,b) c_\lambda(T)
  \end{equation}
and
  \begin{equation} \label{subt-from-XT}
  \subt_T(q,r) ~=~ \sum_{i=1}^{n-1} \sum_{j=1}^i q^i r^j
  \sum_{\lambda\partn n} \phi(\lambda,i,j) c_\lambda(T).
  \end{equation}
\end{thm}

It follows that the chromatic symmetric function is at least as strong an invariant as the
subtree and connector polynomials.  In particular, the path and degree sequences of~$T$
can be recovered from $\CSF_T$, as announced previously in~\cite{Marwag}; this generalizes
a earlier result of Fougere~\cite[Theorem~3.3.1]{Fougere}.  Section~\ref{main-thm-section}
contains the proof of Theorem~\ref{main-thm}, as well as explicit formulas
for the path and degree sequences, and a reinterpretation of \eqref{conn-from-XT} and
\eqref{subt-from-XT} in terms of the usual scalar product on the space of
symmetric functions.

\begin{figure}
  \begin{center}
  \resizebox{3.3in}{0.6in}{\includegraphics{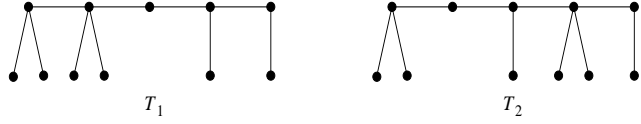}}
  \end{center}
  \caption{\label{counterexample-figure}
           Two trees with the same subtree polynomial
           but different chromatic symmetric functions.}
  \end{figure}

Theorem~\ref{main-thm} implies that $\CSF_T$ is a stronger invariant than $\subt_T$.
In fact, it is \emph{strictly} stronger: the two trees shown in
Figure~\ref{counterexample-figure} have different chromatic symmetric
functions, but the same subtree polynomial.  (Eisenstat and Gordon \cite{EG} constructed an infinite family of pairs of
  non-isomorphic trees with the same subtree polynomials, of which
  Figure~\ref{counterexample-figure} is the smallest example.)
Thus Stanley's question remains open.

As another application of the combinatorial interpretation of the coefficients
$c_\lambda(T)$, we identify some classes of trees
for which the chromatic symmetric function is in fact a complete invariant.
These trees include all \emph{spiders} (trees having exactly one vertex of degree 
$\geq 3$) and some \emph{caterpillars} (trees from which deleting all leaves yields a path,
such as those in Figure~\ref{counterexample-figure}).
We prove in Section~\ref{spider-section} that every spider can be reconstructed from its 
subtree
polynomial, hence from its chromatic symmetric function (generalizing results
of Fougere \cite{Fougere}).
The corresponding problem
for caterpillars is more difficult; however, certain special kinds
of caterpillars can indeed be reconstructed from their chromatic symmetric functions,
and the methods we use to prove this may be extendible to all caterpillars.
%To our knowledge, these are the first known infinite families 
%of trees (other than the trivial instances of paths and stars)
%for which such a reconstruction is known to be possible.

A \emph{unicyclic graph} is a graph with one cycle.  Connected unicyclic graphs
can be recognized as such from their chromatic symmetric functions.  While the
combinatorial data provided by Stanley's expansion \eqref{Stanley's-expansion} is not
as fine for unicyclic graphs as it is for trees, we can still obtain some
uniqueness results for special unicyclic graphs by mimicking our results for
spiders and caterpillars.  In particular, we show in Section~\ref{squid-section}
that no two \emph{squids} (unicyclic graphs with at most one vertex
of degree two or more) can have the same chromatic symmetric function, although
it is not clear whether membership in the class of squids can be determined
from $\CSF_G$.  An analogous result holds for \emph{crabs} (unicyclic
graphs in which every vertex not on the cycle is a leaf) satisfying an additional
technical condition.

\subsection*{Acknowledgements}
Our collaboration began at the Graduate Student
Combinatorics Conference held at the
University of Minnesota on April 16 and 17, 2005.  We thank the
organizers of the conference for their efforts, and we thank
Fran\c{c}ois~Bergeron, Tom~Enkosky, Gary~Gordon, Brandon~Humpert, Rosa~Orellana, Victor~Reiner,
Bruce~Sagan, and Stephanie~van~Willigenburg for
many helpful conversations.  The database of trees constructed by
Piec, Malarz, and Kulakowski \cite{PMK}, and
John~Stembridge's freely available Maple package {\tt SF},
were invaluable for calculating examples and formulating
(and checking small cases of) Conjectures~\ref{positivity-conjecture}
and~\ref{integrality-conjecture}.  
Finally, we thank two anonymous referees
for their suggestions, and in particular for making us aware
of Fougere's undergraduate thesis \cite{Fougere}.

%=================================================================
\section{Background}
\label{background-section}

We assume that the reader is familiar with basic facts about graphs
and trees (see, e.g., \cite[Chapter~I]{Bollobas}).  We denote
a graph $G$ by an ordered pair $(V,E)$, where $V=V(G)$ is the set
of vertices and $E=E(G)$ is the set of edges.  All our graphs
are \emph{simple}; that is, we forbid loops and parallel edges.
The \emph{order} of a graph is its number of vertices.  
A \emph{tree} is a graph $G$ which is acyclic and connected and for which
$\#V(G)=\#E(G)+1$; any two of these conditions together imply the third.
We consider the graph with one vertex and no edges to be a tree, the \emph{trivial tree};
unless otherwise specified, all our statements about trees include this possibility.
A \emph{leaf} of a tree is a vertex of degree~1, that is, with exactly
one incident edge.  Every nontrivial tree has at
least two leaves \cite[p.~11]{Bollobas}.
It is often notationally convenient to ignore the distinction
between a graph and its edge set.

We now review some facts about symmetric functions
(for which the standard references are \cite{Macdonald} and \cite[Chapter~7]{ECII})
and about the chromatic symmetric function (introduced by Stanley in~\cite{Stanley}).

A \emph{partition} is a sequence $\lambda=(\lambda_1,\dots,\lambda_\ell)$
of positive integers in weakly decreasing order.   The numbers $\lambda_k$
are called the \emph{parts} of $\lambda$.  We say that $\lambda$
is a partition of $n$, written $\lambda\partn n$, if $\sum_k \lambda_k=n$.
The number $\ell=\ell(\lambda)$ is called the \emph{length} of $\lambda$.

Let $x_1,x_2,\dots$ be a countably infinite set of commuting indeterminates.
For $k\in\Pp$, the \emph{$k^{th}$ power-sum symmetric function} is
  $$p_k = \sum_{i\geq 1} x_i^k$$
and for a partition $\lambda$ we define
  $$
  p_\lambda = \prod_{k=1}^{\ell(\lambda)} p_{\lambda_k}.
  $$
It is well known that $\{p_\lambda \st \lambda\partn n\}$ is a
basis for the $\Qq$-vector space $\Lambda_n$ consisting of all symmetric functions that
are homogeneous of degree~$n$.

\begin{figure}
\begin{center}
\resizebox{1.6in}{0.4in}{\includegraphics{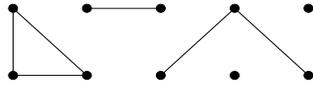}}
\end{center}
\caption{\label{type-figure}
  An edge set of type $(3,3,2,1,1)$.}
\end{figure}

Let $G$ be a graph of order~$n$.
Stanley \cite[Thm.~2.5]{Stanley} proved that
  \begin{equation} \label{Stanley's-expansion}
  \CSF_G = \sum_{A\subseteq E(G)} (-1)^{\#A} p_{\type(A)},
  \end{equation}
where $\type(A)$ is the partition whose parts
are the orders of the connected components of the subgraph of $G$ induced by~$A$
(see Figure~\ref{type-figure} for an example).  Note that $\type(A)$ depends
upon $A$ and $V(G)$, but not on $E(G)$.
We write $c_\lambda(G)$ for the coefficient of $p_\lambda$ in
the expansion \eqref{Stanley's-expansion}; that is,
  \begin{equation} \label{define-c-lambda}
  \CSF_G = \sum_{\lambda\partn n} c_\lambda(G) p_\lambda.
  \end{equation}
We will abbreviate $c_\lambda=c_\lambda(G)$ when no confusion can arise.

The chromatic symmetric function of a graph provides much more
combinatorial information when that graph is a tree.
In general, the coefficient $c_\lambda$
does not count edge sets of type $\lambda$, because $(-1)^{\#A}$ is
not constant for all such sets~$A$.  On the other hand, if $T=(V,E)$ is a tree of order~$n$,
then every $A\subseteq E$ is acyclic, so its induced subgraph has $n-\#A$ connected
components.  Hence $\ell(\type(A))=n-\#A$, and we obtain a useful
combinatorial interpretation for the numbers $c_\lambda$:
  \begin{equation} \label{interpret-c}
  c_\lambda = (-1)^{n-\ell(\lambda)} \#\{A\subseteq E \st \type(A)=\lambda \}.
  \end{equation}
The invariants $c_\lambda$ are far from independent; in particular,
\eqref{interpret-c} implies that
  \begin{equation} \label{given-size}
  \sum_{\lambda:\ \ell(\lambda)=k} (-1)^{n-\ell(\lambda)} c_\lambda(T) = \binom{n-1}{k}.
  \end{equation}

We next list some basic invariants of graphs that can be recovered
from its chromatic symmetric function.  Several of these facts
were previously observed by the second author in \cite{Morin}.
For notational simplicity, we shall often omit
the parentheses and singleton parts of a partition, for instance,
writing $c_j$ rather than $c_{(j,1,1,\dots,1)}$.

\begin{prop} \label{basic-facts}
For every graph $G=(V,E)$,
\begin{itemize}
\item[(i)] the symmetric function $\CSF_G$ is homogeneous of degree $\#V$;
\item[(ii)] $-c_2(G) = \#E$; and
\item[(iii)] the number $k$ of connected components of $G$ is
$\min\{\ell(\lambda) \st c_\lambda(G)\neq 0\}$.
\end{itemize}
If $G$ is a tree, then in addition 
\begin{itemize}
\item[(iv)] for $k\geq 2$, $|c_j(G)|$ is the number of subtrees of $G$ with $j$ vertices; and
\item[(v)] $c_{n-1}(G)$ is the number of leaves of $G$.
\end{itemize}
\end{prop}

\begin{proof}
Assertion (i) is immediate from the definition of $\CSF_G$, and
(ii) follows easily from \eqref{Stanley's-expansion}.  For (iii),
\eqref{Stanley's-expansion} implies that $c_\lambda=0$ whenever
$\ell(\lambda)<k$.  On the other hand,
  $c_\lambda = \sum_A (-1)^{\#A} = (-1)^{n-k} t_G(1,0)$,
the sum over all $A\subseteq E(G)$ with $\lambda(A)=\lambda(E(G))$, where
$t_G$ is the Tutte polynomial of $G$ (see \cite[Ch.~X]{Bollobas}).
Up to sign, this formula counts acyclic orientations of $G$ with exactly
one sink in each component (this follows from \cite[Thm.~8, p.~348]{Bollobas}).
In particular, $c_\lambda\neq 0$.

Assertion (iv) holds because $A\subseteq E(G)$ is
the set of edges of a $j$-vertex subtree if and only if $\type(A)=(j,1,\dots,1)$,
and (v) follows because every subtree of order $n-1$ is of the form
$G-v$, where $v$ is a leaf.
\end{proof}

By (i), (ii) and (iii) of Proposition~\ref{basic-facts},
trees can be distinguished from non-trees by their chromatic symmetric functions.
Moreover, part (v) implies
that paths (trees with exactly two leaves) and stars (trees with exactly one
non-leaf) are determined up to isomorphism by their chromatic symmetric functions.

The \emph{girth} of a graph $G$ is defined as the length of the smallest cycle in $G$,
or $\infty$ if $G$ is acyclic.  With a little more work, we can compute the
girth of $G$ from $\CSF_G$.  The idea is to find the smallest edge
set for which \eqref{given-size} fails.

\begin{prop} \label{girth-prop}
Let $G=(V,E)$ be a graph with $n$ vertices and $m$ edges.  Let
$k$ be the largest number such that
  $\sum_{\lambda\partn n,\:\ell(\lambda)=k} c_\lambda(G)
  ~\neq~ (-1)^{n-k} \binom{m}{n-k}.$
Then the girth of $G$ is $n-k+1$.
\end{prop}

\begin{proof}
Let $g$ be the girth of $G$.
Suppose first that $k>n-g+1$.  Then $n-k<g-1$, so every subset $A\subseteq E$
with $n-k$ edges is acyclic and hence has $k$ connected components.
On the other hand, if $\ell(\type(A))=k$,
then the maximum size of a component of $A$ is $n-(k-1)<g$, so $A$ must be acyclic
and hence must have $n-k$ edges.  Therefore
  $$
  \sum_{\substack{\lambda\partn n\\ \ell(\lambda)=k}} c_\lambda(G)
  ~=~ \sum_{\substack{A\subseteq E\\ \ell(\type(A))=k}} (-1)^{\#A}
  ~=~ \sum_{\substack{A\subseteq E\\ \#A=n-k}} (-1)^{\#A}
  ~=~ (-1)^{n-k}\binom{m}{n-k}.
  $$

Now suppose that $k=n-g+1$.  We claim that $A\subseteq E$ has $k$ components
if and only if it either has $n-k$ edges (hence is acyclic) or is precisely a cycle
of length~$g$.  The ``if'' direction is evident.  For the ``only if'' direction,
suppose that $A$ has $k$ components and is not acyclic, hence contains a cycle $C$.
By definition of $G$, the length of $C$ cannot be less than $g$; on the other hand,
there are at least $k-1$ vertices that do not belong to $C$ (one for each other component
of $A$), so $\#V(C)\leq n-(k-1)=g$.  Thus $C$ has length exactly~$g$.
Moreover, $A-C$ cannot contain any other edge with an
endpoint outside $C$ (because then it would have fewer than $k$ components)
or an edge joining two vertices of~$C$ (because then $G$ would contain a cycle
of length $<g$).  Hence $A=C$ as desired.  Denoting by $\Gamma$ the set
of $g$-cycles of $G$, we have
  \begin{align*}
  \sum_{\substack{\lambda\partn n\\ \ell(\lambda)=k}} c_\lambda(G)
  ~&= \sum_{\substack{A\subseteq E\\ \#A=n-k}} (-1)^{\#A}
    + \sum_{A\in\Gamma} (-1)^{\#A}\\
  &= (-1)^{n-k}\binom{m}{n-k} + (-1)^{n-k+1} \#\Gamma
  ~\neq~ (-1)^{n-k}\binom{m}{n-k}
  \end{align*}
as desired.
\end{proof}

%=================================================================
\section{Proof of the main theorem}
\label{main-thm-section}
%=================================================================

Theorem~\ref{main-thm} expresses the subtree polynomial $\subt_T$
and connector polynomial $\conn_T$ of a tree~$T$ in terms of the
chromatic symmetric function $\CSF_T$.
The first step is to show that $\subt_T$ and $\conn_T$ are interchangeable.
In what follows, we will often abuse notation by ignoring
the distinction between a tree~$T=(V,E)$ and its edge set $E$.

\begin{prop} \label{subt-conn-equivalence}
Let $T$ be a tree.  Then:
\begin{enumerate}
\item $\subt_T(q,r) = \conn_T(qr,q(1-r))$.
\item $\conn_T(x,y) = \subt_T(x+y,x/(x+y))$.
\end{enumerate}
\end{prop}

\begin{proof}
For each nontrivial subtree $S\subset T$, write $L(S)$ for the set of leaf edges of $S$.
Note that $\#L(S)\geq 1$, with equality if and only if $S$ consists of a single edge.
Moreover, observe that $A\cup K(A)=S$ if and only if $L(S)\subseteq A\subseteq S$.
Hence
\begin{eqnarray*}
  \conn_T(qr,q(1-r)) &=& \sum_{A\subseteq T} (qr)^{\#A} (q(1-r))^{\#K(A)}\\
  &=& \sum_{A\subseteq T} q^{\#(A\cup K(A))} r^{\#A} (1-r)^{\#K(A)}\\
  &=& \sum_{\text{subtrees } S\subseteq T}
         q^{\#S} \sum_{\substack{A:\\ L(S)\subseteq A\subseteq S}} r^{\#A} (1-r)^{\#S-\#A}\\
  &=& \sum_S q^{\#S} r^{\#L(S)} \sum_{G\subseteq S-L(S)} r^{\#G} (1-r)^{\#(S-L(S))-\#G}\\
  &=& \sum_S q^{\#S} r^{\#L(S)} \left( r+(1-r) \right)^{\#(S-L(S))} ~=~ \subt_T(q,r),
\end{eqnarray*}
giving the first equality.  Meanwhile,
solving the equations $x=qr,y=q(1-r)$ for $q$ and $r$
yields $q=x+y$, $r=x/(x+y)$, giving the second equality.
\end{proof}

We now prove the main theorem.  To do so, we establish a formula
for the connector polynomial of a tree in terms of its chromatic symmetric function, then
apply Proposition~\ref{subt-conn-equivalence} to obtain a formula for the subtree polynomial.

\begin{proof}[Proof of Theorem~\ref{main-thm}.]
By definition, the coefficient of $x^ay^b$ in $\conn_T(x,y)$ is
  \begin{multline*}
  \#\{A\subseteq T \st \#A=a,\; \#K(A)=b\}
  ~=~ \sum_{\substack{A\subseteq T\\ \#A=a,\; \#K(A)=b}} 1\\
   =~ (-1)^b \sum_{\substack{A\subseteq T\\ \#A=a}}\ \ 
  \sum_{\substack{B\subseteq T-A\\ \#B=b}}\ \ 
  \sum_{\substack{C\subseteq B\\ K(A)\subseteq C}} (-1)^{\#C},
  \end{multline*}
because the innermost sum vanishes unless $B=K(A)$, when it is $(-1)^b$.  Setting $D=B-C$,
we may rewrite this expression as
  $$
  (-1)^b \sum_{\substack{A\subseteq T\\ \#A=a}}\ \ 
  \sum_{\substack{C\subseteq T-A\\ K(A)\subseteq C}}\ \ 
  \sum_{\substack{D\subseteq T-A-C\\ \#D=b-\#C}} (-1)^{\#C}
  $$
and setting $F=C\cup A$ yields
\begin{align}
  (-1)^b &
  \sum_{\substack{A\subseteq T\\ \#A=a}}\ \ 
  \sum_{\substack{F\subseteq T\\ A\cup K(A)\subseteq F}}\ \ 
  \sum_{\substack{D\subseteq T-F\\ \#D=a+b-\#F}} (-1)^{\#F-a} \notag\\
  &=~
  (-1)^{a+b} \sum_{F\subseteq T}\ 
  \sum_{\substack{A:\ \#A=a,\\ A\cup K(A)\subseteq F}}
  \binom{\#(T-F)}{a+b-\#F} (-1)^{\#F} \notag\\
  &=~
  (-1)^{a+b} \sum_{\lambda\partn n}\ 
  \sum_{\substack{F\subseteq T\\ \type(F)=\lambda}}
  \binom{\#(T-F)}{a+b-\#F}
  \sum_{\substack{A:\ \#A=a,\\ A\cup K(A)\subseteq F}} (-1)^{\#F} \notag\\
  &=~
  (-1)^{a+b}
  \sum_{\lambda\partn n}
  \binom{\ell(\lambda)-1}{\ell(\lambda)-n+a+b}
  (-1)^{n-\ell(\lambda)}
  \sum_{\substack{F\subseteq T\\ \type(F)=\lambda}} \alpha(F) \label{almost}
\end{align}
where
  $\alpha(F) = \#\{A \st \#A=a,\: A\cup K(A)\subseteq F\}.$
The set $A\cup K(A)$ is connected, so if it is a subset of $F$ then it must be a subset of some component
of $F$.  On the other hand, if $F'$ is a (possibly trivial) component of $F$ and $A\subseteq F'$, then
$A\cup K(A)\subseteq F'$, because $F'$ is a tree containing $A$ and $A\cup K(A)$ is the unique
minimal such tree.  Thus if $\type(F)=\lambda$ then
  \begin{equation} \label{alpha-formula}
  \alpha(F) = \sum_{k=1}^{\ell(\lambda)} \binom{\lambda_k-1}{a}.
  \end{equation}
Note that this formula is valid only if $a>0$.
Substituting \eqref{alpha-formula} into \eqref{almost}, we obtain
  \begin{align*}
  & (-1)^{a+b}
  \sum_{\lambda\partn n}
  \binom{\ell(\lambda)-1}{\ell(\lambda)-n+a+b}
  (-1)^{n-\ell(\lambda)}
  \sum_{\substack{F\subseteq T\\ \type(F)=\lambda}}
  \sum_{k=1}^{\ell(\lambda)} \binom{\lambda_k-1}{a} \\
  &=~
  (-1)^{a+b}
  \sum_{\lambda\partn n}
  \binom{\ell(\lambda)-1}{\ell(\lambda)-n+a+b}
  \sum_{k=1}^{\ell(\lambda)} \binom{\lambda_k-1}{a}
  c_\lambda(T) \\
  &=~
  \sum_{\lambda\partn n} \psi(\lambda,a,b)  c_\lambda(T)
  \end{align*}
(where $\psi(\lambda,a,b)$ is defined by \eqref{define-psi}),
giving the desired formula \eqref{conn-from-XT}.

\smallskip

We now turn to the proof of \eqref{subt-from-XT}.
By Proposition~\ref{subt-conn-equivalence}, we have
\begin{multline*}
\subt_T(q,r) \;=\; \conn_T(qr,q(1-r))\\
=\;
\sum_{a>0} \sum_{b\geq 0} (qr)^a (q(1-r))^b (-1)^{a+b}
\sum_{\lambda\partn n} \binom{\ell(\lambda)-1}{\ell(\lambda)-n+a+b}
  \sum_{k=1}^{\ell(\lambda)} \binom{\lambda_k-1}{a} c_\lambda(T).
\end{multline*}
Setting $i=a+b$, we may rewrite the last expression as
$$
\sum_{i>0}\sum_{a=1}^i (-1)^i q^i r^a (1-r)^{i-a}
\sum_{\lambda\partn n} \binom{\ell(\lambda)-1}{\ell(\lambda)-n+i}
  \sum_{k=1}^{\ell(\lambda)} \binom{\lambda_k-1}{a} c_\lambda(T).
$$
Applying the binomial expansion to $(1-r)^{i-a}$ yields
$$
\sum_{i>0}\sum_{a=1}^i (-1)^i q^i r^a
\sum_{h=0}^{i-a} \binom{i-a}{h} (-1)^h r^h
\sum_{\lambda\partn n} \binom{\ell(\lambda)-1}{\ell(\lambda)-n+i}
\sum_{k=1}^{\ell(\lambda)} \binom{\lambda_k-1}{a} c_\lambda(T).
$$
Now setting $h=j-a$ gives
$$
\sum_{i>0}\sum_{a=1}^i (-1)^i q^i
\sum_{j=a}^i \binom{i-a}{j-a} (-1)^{j-a} r^j
\sum_{\lambda\partn n} \binom{\ell(\lambda)-1}{\ell(\lambda)-n+i}
\sum_{k=1}^{\ell(\lambda)} \binom{\lambda_k-1}{a} c_\lambda(T)
$$
and setting $a=d$ and rearranging gives
\begin{multline*}
\sum_{i>0}\sum_{j=1}^i q^i r^j
\sum_{\substack{\lambda\partn n\\ \ell(\lambda)=\ell}}
\left(
  (-1)^{i+j} \binom{\ell-1}{\ell-n+i}
  \sum_{d=1}^j (-1)^d \binom{i-d}{j-d}
  \sum_{k=1}^{\ell} \binom{\lambda_k-1}{d}
\right) c_\lambda(T)\\
=~
\sum_{i>0}\sum_{j=1}^i q^i r^j
\sum_{\substack{\lambda\partn n\\ \ell(\lambda)=\ell}}
\phi(\lambda,i,j) c_\lambda(T)
\end{multline*}
which is the desired formula \eqref{subt-from-XT}.
\end{proof}

Two basic invariants of a tree are its \emph{path sequence} and its
\emph{degree sequence}.  The path sequence of $T$ is defined as
$(\pi_1,\pi_2,\dots)$, where $\pi_i=\pi_i(T)$ is the number of $i$-edge paths in $T$.
The degree sequence of $T$ is defined as $(\delta_1,\delta_2,\dots)$,
where $\delta_j=\delta_j(T)$ is the number of degree-$j$ vertices in $T$.
Knowing the degree sequence is equivalent to knowing the \emph{star sequence}
$(\sigma_1,\sigma_2,\dots)$, where $\sigma_k=\sigma_k(T)$ is the number
of $k$-edge stars in $T$.  Indeed, it is not hard to see that
  $$\sigma_k = \sum_{j\geq k} \binom{j}{k} \delta_j$$
for every $2\leq k\leq n-1$, and so
  $$\delta_j = \sum_{k\geq j} \binom{k}{j} (-1)^{j+k} \sigma_k.$$

\begin{cor} \label{degree-and-path-sequences}
The degree and path sequences of a tree $T$ can be recovered from
its chromatic symmetric function.
\end{cor}

\begin{proof}
The key observation, due to Chaudhary and Gordon \cite[Proposition~18]{CG},
is that the path and star sequences of~$T$ can be recovered from $\subt_T$.
Indeed, $\pi_1$ is the number of edges of~$T$, and
for every $i\geq 2$, $\pi_i$ is just the coefficient of $q^ir^2$ in $\subt_T(q,r)$.
Meanwhile, for every $k\geq 1$, $\sigma_k$ is the coefficient of $q^kr^k$.
\end{proof}

We note that Fougere had proved~\cite[Theorem~3.3.1]{Fougere} that the
sum of the squared vertex degrees, $\sum_j \delta_jj^2$, could be obtained
from the coefficient of the monomial symmetric function $m_{(3,1,1,\dots)}$
in $\CSF_T$.

We can rephrase the formulas for $\conn_T$ and $\subt_T$ in terms of the usual scalar product
$\langle\cdot,\cdot\rangle$ on the space $\Lambda_n$ of degree-$n$ symmetric functions
(see \cite[\S7.9]{ECII} or \cite[\S I.4]{Macdonald}), where $n$ is the order of~$T$.
Define symmetric functions $\Psi_n(x,y)$ and $\Phi_n(q,r)$ by
  \begin{align*}
  \Psi_n(x,y) &= \sum_{a>0} \sum_{b\geq 0} x^a y^b
  \sum_{\lambda\partn n} \psi(\lambda,a,b) \frac{p_\lambda}{z_\lambda},\\
  \Phi_n(q,r) &= \sum_{i=1}^{n-1} \sum_{j=1}^i q^i r^j
  \sum_{\lambda\partn n} \phi(\lambda,i,j) \frac{p_\lambda}{z_\lambda}.
  \end{align*}
Then the formulas \eqref{conn-from-XT} and \eqref{subt-from-XT} are
respectively equivalent to
  \begin{align*}
  \conn_T(x,y) &= \scp{\Psi_n(x,y)}{\CSF_T},\\
  \subt_T(q,r) &= \scp{\Phi_n(q,r)}{\CSF_T}.
  \end{align*}

The symmetric function $\Psi_n$ appears to have certain positivity
and integrality properties,
as we now explain.  The \emph{$i^{th}$ homogeneous symmetric function} $h_i$
in indeterminates $x_1,x_2,\dots$ is the sum of all monomials of degree~$i$,
and for a partition $\lambda=(\lambda_1,\dots,\lambda_\ell)$ we define
$h_\lambda=h_{\lambda_1}\cdots h_{\lambda_\ell}$.  The symmetric functions
$\{h_\lambda \st \lambda\partn n\}$ form a vector space basis for $\Lambda_n$
\cite{Macdonald}, so there is a unique list of rational numbers
$\xi(\lambda,i,j)\in\Qq$ such that
  $$\Psi_n(x,y) = \sum_{i,j} \sum_{\lambda\partn n} \xi(\lambda,i,j) x^i y^j h_\lambda.$$

\begin{conj}[Positivity] \label{positivity-conjecture}
Let $\mu\partn n$ be a partition, and let
$\epsilon(\mu)$ be the number of parts of $\mu$ of even length.
Then, for all integers $i,j$,
  $(-1)^{\epsilon(\mu)} \xi(\mu,i,j) \geq 0.$
\end{conj}

\begin{conj}[$z$-Integrality] \label{integrality-conjecture}
Let $\mu\partn n$ be a partition.  Then, for all integers $i,j$,
the number $\xi(\mu,i,j)z_\mu$ is an integer.
\end{conj}

We have verified
Conjectures~\ref{positivity-conjecture} and~\ref{integrality-conjecture}
computationally\footnote{
  A Maple worksheet containing the calculations is available at\\
  \url{http://math.ku.edu/~jmartin/sourcecode/}.}
for all $n\leq 20$, which we think is strong evidence that they hold for all~$n$.
A formula for $\xi(\mu,i,j)$ can be written out explicitly using the known
transition matrices between bases of symmetric functions (see~\cite{ER}).
However, we do not know a direct combinatorial interpretation for
$\xi(\mu,i,j)$ or for $\xi(\mu,i,j)z_\mu$.

One might hope that for every two trees $T,U$ with $\#V(T)>\#V(U)>1$, the number of subtrees of~$T$
isomorphic to~$U$ might be given by a scalar product $\scp{\zeta_U}{\CSF_T}$,
where $\zeta_U$ is some symmetric function independent of~$T$.  Such a result
would generalize Corollary~\ref{degree-and-path-sequences} (which covers only
the case that~$U$ is a path or a star) and, by a theorem of Harary and Palmer
\cite{HP}, would imply that every tree is distinguished by its chromatic
symmetric function.  In fact, it appears that such a function $\zeta_U$ exists
\emph{only if} $U$ is a star or a path, as we have verified computationally
for all $U$ of order~$\leq 8$, with one trivial exception.\footnote{
  Up to isomorphism, there are three four-edge trees: the star $S_4$,
  the path $P_4$, and another tree~$U$.  Since the number of four-edge
  subtrees of~$T$ is just $c_5(T)$,  we have $\zeta_U=c_5-\zeta_{S_4}-\zeta_{P_4}$.}

Theorem~\ref{main-thm} does not resolve Stanley's question, because $\subt_T$
is not a complete isomorphism invariant.  Indeed, the two
trees $T_1,T_2$ shown in Figure~\ref{counterexample-figure} share
the same subtree polynomial; this is a special case of a theorem of
Eisenstat and Gordon~\cite{EG}.
On the other hand, $\CSF_{T_1}\neq\CSF_{T_2}$.
This inequality follows from Tan's calculations \cite{Tan},
and also for the following elementary reason.
Let $A\subset E(T_1)$ be the edge set obtained by deleting the two
rightmost horizontal edges in Figure~\ref{counterexample-figure};
then $\type(A)=(7,2,2)$.  On the other hand, no subset of
$E(T_2)$ has that type.  Therefore
$c_{(7,2,2)}(T_1)\neq 0$ and $c_{(7,2,2)}(T_2)=0$.

The remainder of the article is devoted to identifying special classes of trees~$T$
for which the invariants $c_\lambda$ suffice to reconstruct~$T$ up to isomorphism.

%=================================================================
\section{The Chromatic Symmetric Function Distinguishes Spiders}
\label{spider-section}
%=================================================================

A tree is a \defterm{spider} (or \emph{starlike tree})
if exactly one of its vertices has degree $\geq 3$.
By Corollary~\ref{degree-and-path-sequences}, whether or not a tree is
a spider can be determined from its subtree polynomial.
A spider may equivalently be defined as a collection of edge-disjoint paths
(the \defterm{legs}) joined at a common endpoint $t$ (the \defterm{torso}).

Up to isomorphism, every spider on~$n$ vertices can be described
by a partition $\lambda\partn n-1$ whose parts are the lengths of its legs (so
$\ell(\lambda)\geq 3$).  We denote the corresponding
spider by $T_\lambda$.  Note that $\ell(\lambda)$ equals both the number
of leaves of $T_\lambda$ and the degree of its torso.

We will show that the isomorphism type of a spider can be determined
from its subtree polynomial, hence from its chromatic symmetric function.
Fougere \cite[Chapter~2]{Fougere} had previously shown that \emph{forks} (spiders
with exactly one leg of length $>1$ and \emph{extended stars} (spiders
in which every leg has length $k$ or $k+1$ for some $k$) could be
reconstructed from their chromatic symmetric functions.

Before continuing, we describe a combinatorial
problem whose solution will play a role in the proof.
Let $m_1,\dots,m_k$ be nonnegative integers with $\sum m_i = \ell$.
Suppose that we have $\ell$ distinguishable boxes, of which $m_i$ have
capacity~$i$ for each $i\in[k]$.  Let $\Omega(m_1,\dots,m_k)$ be
the number of ways of distributing $k$ indistinguishable
balls among these boxes so that no box is filled beyond its capacity.
In general, it is
not easy to write down a closed formula for $\Omega(m_1,\dots,m_k)$,
although individual instances can be computed using an inclusion-exclusion
argument (for example).

\begin{thm} \label{spider-thm:1}
Let $\lambda\partn n-1$ be a partition with $\ell=\ell(\lambda)\geq 3$, and
let $T=T_\lambda$ be the corresponding spider.  Then $T$ can be reconstructed
from its subtree polynomial.
\end{thm}

\begin{proof}
First, suppose that $\ell=3$.  Then $\lambda_1+\lambda_2$, $\lambda_1+\lambda_2+\lambda_3$,
and $\lambda_1\lambda_2\lambda_3$ are respectively the diameter, number of edges, and
number of three-leaf subtrees of $T_\lambda$.  These invariants can be recovered from
$\subt_{\lambda}$, and together they determine $\lambda$.

Now suppose that $\ell>3$.  Let $m_k$ denote the
number of parts of size $k$ in $\lambda$,
and let $s(i,j)$ denote the number of subtrees of the spider
$T$ with $i$ edges and $j$ leaf edges (that is, the coefficient
of $q^ir^j$ in $\subt_T(q,r)$).  We will show by induction on~$k$ that
$m_k$ can be calculated from the numbers $s(i,j)$.

First, suppose $k=1$.
Since there is a bijection between legs of $T$ of
length 1 (i.e., consisting of a single edge) and subtrees with $n-2$ edges and $\ell-1$ legs
(which are formed by deleting such an edge).  Hence $m_1=s(n-2,\ell-1)$.

For $k>1$, we can choose a subtree $S\subset T$ with $n-1-k$ edges and $\ell-1$ leaves
as follows.  First, fix $j\in[k]$ and delete a leg of length $j$;
there are $m_j$ ways to do this.  If $j<k$, then we still need to delete $k-j$ more edges.
It suffices
to specify how many edges to delete from the end of each of the other $\ell-1$ legs, so
the number of ways to do this is
the solution to the balls-in-boxes problem described above, regarding
the $k-j$ edges to be deleted as balls and each remaining leg of length $i$ as
a box of capacity $\min(i-1,k-j)$ (since deleting the entire leg will result in a tree
with fewer than $\ell-1$ leaves).  Therefore  $s(n-k-1,\ell-1)$ is given by the formula
\begin{multline} \label{omega-eqn}
 m_k + \sum_{j=1}^{\lfloor k/2 \rfloor} m_j
     \Omega\big(m_2,\dots,m_{j-1},m_j-1,m_{j+1},\dots,m_{k-j},
     \ell-(m_1+\cdots+m_{k-j})\big) \\
  + \sum_{j=\lfloor k/2 \rfloor+1}^{k-1} m_j
     \Omega\big(m_2,\;\dots,\;m_{k-j},\;\ell-1-(m_1+\cdots+m_{k-j})\big).
  \end{multline}
By induction, the $s(i,j)$ determine $m_1,\dots,m_{k-1}$, and
equation \eqref{omega-eqn} implies that they determine $m_k$ as well.
\end{proof}

Another way of reconstructing a spider from its chromatic symmetric function
will be useful in Section~\ref{squid-section}.
Define the \emph{two-part portion} $ \TPCSF_G$ of $\CSF_G$ by
  \begin{equation} \label{two-part}
  \TPCSF_G = \sum_{\ell(\lambda)=2} c_\lambda p_\lambda.
  \end{equation}
While $\TPCSF_G$ is evidently a much weaker invariant than $\CSF_G$, it contains enough
information to distinguish among spiders.

\begin{thm} \label{spider-thm:2}
Let $T$ be a spider of order $n$.  Let $m=\lfloor\frac{n}{2}\rfloor$,
let $\ep=n-2m$, 
let $d_a=|c_{a,n-a}(T)|$ for $1\leq a \leq m$, and let $d(T)=(d_1,\dots,d_m)$.
Then one of the following conditions holds:
\begin{enumerate}
\item The sequence $d(T)$ is a partition,
that is, $d_1\geq d_2\geq\cdots\geq d_m\geq 0$.
\item There is a number $t\leq m$ such that $d_1\geq\cdots\geq d_t=1$.
Moreover, $d_{t+1}=d_{t+2}=\cdots=d_{m-1}=2$, and $d_m\in\{1,2\}$.
\end{enumerate}
In the first case, let $\mu=d(T)$.  In the second case,
define a partition $\mu$ from $d(T)$ by replacing all 2's after the $t^{th}$ place with
two 1's.

Then $T=T_\lambda$, where $\lambda$ is the conjugate partition to $\mu$.
\end{thm}

\newcommand{\Clittle}{C_{\text{little}}}
\newcommand{\Cbig}{C_{\text{big}}}

\begin{proof}
For each edge $e\in E(T)$, denote by $\|e\|$ the minimum of the
orders of the components of $T-e$, so that
$d_a=\#\{e\st \|e\|=a\}$.  If we label the legs
of $T$ as $L_1,\dots,L_\ell$, where $L_i$ has $\lambda_i$
edges, and label the edges of each $L_i$ as $e_{i,1},e_{i,2},\dots,
e_{i,\lambda_i}$, starting from the leaf and reading in toward the torso,
then $\|e_{i,j}\|=\min(j,n-j)$.

\Case{1}{$\lambda_1 \leq n-\lambda_1$}
Then $j\leq n-j$ for every edge $e_{i,j}$.  In this case $d_a = \#\{i \st \lambda_i\geq a\}$,
and $d(T)$ is just the conjugate partition of $\lambda$.

\Case{2}{$\lambda_1 > n-\lambda_1$} %%%\lambda_2+\cdots+\lambda_\ell+1$}
%If $n$ is even, then
%  \begin{align*}
%  \|e_{1,1}\|   &= 1,   & \|e_{1,2}\|   &= 2,   & \dots,& & \|e_{1,m}\|         &= m,\\
%  \|e_{1,m+1}\| &= m-1, & \|e_{1,m+2}\| &= m-2, & \dots,& & \|e_{1,\lambda_1}\| &= n-\lambda_1
%  \end{align*}
%while if $n$ is odd then
%  \begin{align*}
%  \|e_{1,1}\|   &= 1, & \|e_{1,2}\|   &= 2,   & \dots,& & \|e_{1,m}\|         &= m,\\
%  \|e_{1,m+1}\| &= m, & \|e_{1,m+2}\| &= m-1, & \dots,& & \|e_{1,\lambda_1}\| &= n-\lambda_1.
%  \end{align*}
If $i\neq 1$, then
$\lambda_i \leq (\lambda_2+\cdots+\lambda_k)-1 = n-\lambda_1-2$,
and $\|e_{i,j}\|=j$ for every $j$.
So we can give a formula for $d(T)$ in terms of $\lambda$:
  \begin{equation} \label{conjugate-complicated}
  d_a = \begin{cases}
    \#\{i \st \lambda_i\geq a\} & \quad\text{if } 1 \leq a \leq n-\lambda_1-1,\\
%    1 & \quad\text{if } a = n-\lambda_1-1,\\
    2 & \quad\text{if } n-\lambda_1 \leq a \leq m-1,\\
    1+\ep & \quad\text{if } a=m.
  \end{cases}
  \end{equation}

Note that $d_{n-\lambda_1-1}=1$, because $e_{1,n-\lambda_1-1}$
is the unique edge whose deletion contains a component of order $n-\lambda_1-1$.
On the other hand, $n-\lambda_1 = 2m+\ep-\lambda_1 \leq 2m+\ep-(m+1) = m+\ep-1$,
so $d_{n-\lambda_1}=2$ whether $n$ is odd or even.  Therefore, $d(T)$ is not
a partition, but has the form described in case~(2) of the theorem.
%  \begin{equation} \label{almost-partition}
%  \begin{aligned}
%  & d_1 \;\geq \;\cdots \;\geq \;d_{n-\lambda_1-1}=1,\;2,\;\dots,\;2,\;1 \qquad\text{or}\\
%  & d_1 \;\geq \;\cdots \;\geq \;d_{n-\lambda_1-1}=1,\;2,\;\dots,\;2.
%  \end{aligned}
%  \end{equation}
Formula~\eqref{conjugate-complicated} implies that if we ``iron out''
$d(T)$ by replacing every 2 after the $(n-\lambda_1)^{th}$ entry with two 1's,
we will obtain the conjugate partition to $\lambda$.

\smallskip
We now see how to recover $\lambda$, and thus the isomorphism type of the
spider $T_\lambda$, from the data $d(T_\lambda)$.  Either $d(T_\lambda)$ is
a partition, in which case it is the conjugate of $\lambda$, or it has the form
just described, in which case the ``ironing-out'' operation
yields the conjugate of $\lambda$.  This is precisely the statement of the theorem.
\end{proof}

%=================================================================
\section{Chromatic Symmetric Functions of Some Caterpillars}
\label{caterpillar-section}
%=================================================================

A \emph{caterpillar} is a tree $T$ with the property that the induced subgraph
on the non-leaf vertices is a nontrivial path, called the \emph{spine} of $T$.  
That is, every vertex of $T$ either lies on the spine, or is a leaf whose unique
neighbor lies on the spine.  Since the spine is nontrivial, a caterpillar must
have at least four vertices.
%  If the vertices of the spine are labeled in order as
%$v_0,v_1,\dots,v_s$, then the ordered list
%
%sequence $d_0=\deg_T(v_0),\dots,d_s=\deg_T(v_s)$, up to reversal,
%determines the isomorphism type of the caterpillar.
%  If this
%sequence is palindromic, that is, $d_i=d_{s-i}$ for all $i$,
%then $T$ is called a \emph{symmetric caterpillar}.

While Eisenstat and Gordon's result in \cite{EG} rules out the
possibility of distinguishing caterpillars by their subtree polynomials,
there is still reason to hope that the additional information provided by the
chromatic symmetric function of a caterpillar may suffice to reconstruct it
up to isomorphism.

Our first result is that the chromatic symmetric function distinguishes
caterpillars from non-caterpillars.  The number of leaves and the diameter
of a tree~$T$ (the maximum length of a path in~$T$)
can be recovered from $\CSF_T$ by Corollary~\ref{degree-and-path-sequences},
so it suffices to prove the following fact.

\begin{prop}
Let $T$ be a tree with~$n\geq 4$ vertices.  Then $T$ is a caterpillar if and only if
$\diam(T)-1 = n-\delta_1(T)$.
\end{prop}

\begin{proof}
If $T$ is a caterpillar, then every path of maximum length consists of the spine
together with a leaf attached to each of its endpoints, hence contains all the
non-leaf vertices and two other (leaf) vertices.  In particular, the number of edges
in such a path is one more than the number of non-leaf vertices.  On the other hand, if
$P$ is a path of maximum length in $T$, then the internal vertices of $P$ are
not leaves of $T$ (because each has two neighbors in $P$) but its endpoints are
(otherwise $P$ could be lengthened).  If $\diam(T)-1 = n-\delta_1(T)$,
then all the vertices not lying on $P$ must be leaves, which is to say that $T$
is a caterpillar.
\end{proof}

Let $T$ be a caterpillar with spine vertices $v_0,\dots,v_s$.  For each~$i$,
let $e_i=\deg(v_i)-1$, where $\deg(v_i)$ denotes the degree of the vertex $v_i$.
Gordon and McDonnell \cite[Lemma~2]{GM} showed that the numbers $e_i$ are almost
determined by the path sequence of $T$, and are indeed determined by the path
sequence when the caterpillar is \emph{symmetric} (that is, $e_i=e_{s-i}$
for $0\leq i\leq s$).  Therefore, every symmetric caterpillar
is determined up to isomorphism by its chromatic symmetric function, a result
proved in another way by the second author \cite[Theorem~4.3.1]{Morin}.

We now describe another class of caterpillars that can be reconstructed
from their chromatic symmetric functions.  We retain the labeling of the
vertices of $T$ as $v_0,\dots,v_s$.  Let $f_i$ be the number of \emph{leaves}
adjacent to $v_i$, so that $f_i=\deg(v_i)-1$ for $i=0,s$ and $f_i=\deg(v_i)-2$
for $0<i<s$.  Call $f_i$ the $i^{th}$ \emph{leaf number} of $T$.  In addition,
call a partition $\lambda$ \emph{singleton-free} if all its parts are
at least 2.

\begin{thm} \label{distinct-caterpillar-thm}
Let $T$ be a caterpillar whose leaf numbers $f_i$ are strictly positive and distinct.
Then $T$ can be reconstructed from its chromatic symmetric function.
\end{thm}

\begin{proof}
Let $L$ be the set of leaf edges of $T$.
Since $f_i>0$ for all $i$, every spine vertex is adjacent to at least one
leaf.  Therefore, the edge sets $A\subseteq T$ such that $\type(A)$ is singleton-free
are precisely those that contain $L$.  In particular, $\lambda=\type(L)$
is the unique singleton-free partition with $s+1$ parts whose coefficient $c_\lambda(T)$
is non-zero.  Up to reordering, the parts of $\lambda$ are the numbers $f_0+1,\dots,f_s+1$.
Furthermore, for each spine edge $e_i=v_{i-1}v_i$,
the edge set $L\cup\{e_i\}$ contributes $(-1)^{n-s}$ to $c_{\mu_i}(T)$, where
   $\mu_i = \{f_0+1,\; \dots,\; f_{i-2}+1,\; f_{i-1}+f_i+2,\; f_{i+1}+1,\;
    \dots,\; f_s+1\}.$
Note that the partitions $\mu_i$ are all distinct.  Moreover, $\mu_i$ has
$s-1$ of its parts in common with $\lambda$; the remaining two parts
of $\lambda$ must be the leaf numbers of adjacent spine vertices of $T$.
(This statement is valid even if the parts of $\mu_i$ are not all distinct.)
In this way we can recover the leaf numbers of the endpoints of every
edge of~$T$, and this data specifies the caterpillar up to isomorphism.
\end{proof}

Using the same argument, we can relax the condition of the theorem slightly:
we need only require that all leaf numbers are positive and that for each~$k$,
the set of spine vertices with leaf number~$k$ form a subpath of the spine.

%=================================================================
\section{Unicyclic Graphs: Squids and Crabs}
\label{squid-section}
%=================================================================

Despite the title of this article, we devote the last section to a family
of graphs $G$ that are not trees, but rather \emph{unicyclic}; that is,
$G$ contains a unique cycle.  Equivalently, $c=n-e+1$, where $c,n,e$ are
the numbers of components, vertices and edges respectively, so
unicyclicity can be detected from $X_G$ by Proposition~\ref{basic-facts}.
While we can no longer interpret the coefficients $c_\lambda$ as in
\eqref{interpret-c}, we can use Stanley's expansion of $\CSF_G$
in terms of broken circuits \cite[Thm.~2.9]{Stanley}.  A special
case of that result is the following:
if $G$ has a unique cycle $C$ and $e_0$ is an edge of $C$, then
  \begin{equation} \label{reducedsum}
  \CSF_G = \sum_{\substack{A \subseteq E(G)\\ C-e_0 \not \subseteq A}}
  (-1)^{\#A} p_{\type(A)}.
  \end{equation}

We do not know whether there exist two unicyclic graphs with the same
chromatic symmetric functions.

A \emph{squid} is a connected unicyclic graph with a single vertex $v$
of degree greater than 2.  Note that $v$ must lie on the cycle.
A squid is described up to isomorphism by the length of
its cycle and of the \emph{tentacles} (the paths
from the leaf vertices to~$v$).  For example, if $T_\lambda$ is the
spider whose leg lengths are given by the parts of $\lambda$,
then adding an edge between the leaves at the ends of the two longest
legs produces a squid with cycle length $\lambda_1+\lambda_2$ and tentacle
lengths $\lambda_3,\lambda_4,\dots$.

It is not clear how to determine from $\CSF_G$ whether or not a unicyclic
graph $G$ is a squid (for instance, we cannot recover the degree sequence
of an arbitrary graph from its chromatic symmetric function as we can for
a tree).  Nevertheless, the following uniqueness result does hold.

\begin{thm} \label{squid-thm}
No two non-isomorphic squids have the same chromatic symmetric function.
\end{thm}

\begin{proof}
Let $S$ be a squid with unique cycle~$C$.  Let $k+1$ be the length of~$C$; this
number can be recovered from $\CSF_S$ by Proposition~\ref{girth-prop}.
Let $v$ be the unique vertex of $S$ of degree $>2$, and let $\mu$ be
the partition whose parts are the edge lengths of the tentacles.
Label the edges of $C$ as $e_0,e_1,\ldots,e_k$, starting at $v$ and proceeding
around the cycle.  By \eqref{reducedsum} and inclusion-exclusion, we obtain
  \begin{align*}
  \CSF_S &= \sum_{i=1}^{k} \left( \sum_{e_i \not\in A} (-1)^{\#A} p_{\type(A)} \right)
     - \sum_{1\leq i<j\leq k} \left( \sum_{e_i,e_j \not\in A} (-1)^{\#A} p_{\type(A)} \right)
     + \dots\\
    &= \sum_{i=1}^{k} \CSF_{S-e_i}
     - \sum_{1 \leq i < j \leq k} \CSF_{S-e_i-e_j}
     + \dots
  \end{align*}
where the omitted terms involve edge sets $A$ lacking three or more edges from $C-e_0$.
Deleting three or more edges from $S$ produces a graph with three or more connected components,
so passing to two-part portions as in \eqref{two-part} yields
  $\TPCSF_S = \sum_{i=1}^{k} \TPCSF_{S-e_i} - \sum_{1\leq i < j \leq k} \TPCSF_{S-e_i-e_j}.$
Each graph $S-e_i-e_j$ has exactly two components, of sizes $j-i$ and $n-j+i$.
Removing additional edges will strictly increase the number of components, so
  $$\TPCSF_S = \sum_{i=1}^{k} \TPCSF_{T_{\mu, i,k-i}}
    - \sum_{1 \leq i < j \leq k} (-1)^{n-2} p_{(j-i, n-j+i)}$$
where $T_{\mu,i,k-i} = S-e_i$ is the spider with legs whose lengths
are $i$, $k-i$, and the parts of $\mu$.  Therefore, $\CSF_S$ determines
the quantity
  \begin{equation} \label{sumofspiders}
  \sum_{i=1}^{k} \TPCSF_{T_{\mu, i,k-i}} =
  \TPCSF_{S} + \sum_{1 \leq i < j \leq k} (-1)^{n-2} p_{(j-i, n-j+i)}.
  \end{equation}

Leaving the foregoing calculations aside for the moment, we note that if $T$ is a
tree with $n-1$ edges, then to calculate $\TPCSF_T$ we need only consider the edge
subsets of cardinality $n-2$.  In particular, if $T=T_\lambda$ is a spider, then
  $\TPCSF_T = (-1)^{n-2} \sum_{i=1}^{\ell(\lambda)} \sum_{j=1}^{\lambda_i} p_{(j,n-j)}.$
It follows that for any partition $\mu$ and numbers $k,i$, we have
  \begin{align*}
  \TPCSF_{T_{\mu,k}} &= \TPCSF_{T_{\mu,i,k-i}} + (-1)^{n-2}\left(
    \sum_{j=1}^{k} p_{(j,n-j)} - \sum_{j=1}^{i} p_{(j,n-j)} - \sum_{j=1}^{k-i} p_{(j,n-j)} \right)\\
  &= \frac{1}{k}\sum_{i=1}^k \left( \TPCSF_{T_{\mu,i,k-i}} + (-1)^{n-2}\left(
    \sum_{j=1}^{k} p_{(j,n-j)} - \sum_{j=1}^{i} p_{(j,n-j)} - \sum_{j=1}^{k-i} p_{(j,n-j)} \right) \right)
  \end{align*}
which can be computed from $\CSF_S$ using \eqref{sumofspiders}.  Meanwhile,
by Theorem~\ref{spider-thm:2}, we can reconstruct the spider $T_{\mu,k}$ from
$\TPCSF_{T_{\mu, k}}$.  In particular, we can reconstruct the partition $\mu$,
which gives the tentacle lengths of the squid $S$.
\end{proof}

Just as squids can be regarded as the unicyclic analogues of spiders,
the unicyclic analogues of caterpillars are \emph{crabs}: connected
unicyclic graphs in which every vertex not lying on the cycle is a leaf.
In analogy to Theorem~\ref{squid-thm} and its proof, we ask whether
is it possible to use the results of Section~\ref{caterpillar-section}
to classify the chromatic symmetric functions of (some) crabs.  As a starting point,
we prove the following analogue of Theorem~\ref{distinct-caterpillar-thm}.

\begin{prop}
Let $G$ be a crab such that the degrees of its non-leaf vertices are all
distinct and all greater than 2.  Then, if $H$ is another crab with this
property, and $\CSF_G=\CSF_H$, then $G\isom H$.
\end{prop}

\begin{proof}
We will show that the coefficients of $\CSF_G$,
together with the knowledge that $G$ is a crab with the property just mentioned,
determine $G$ up to isomorphism.

The girth $g$ of $G$ can be recovered from $\CSF_G$ by
Proposition~\ref{girth-prop}.  Let $C$ be the unique cycle of $G$, and
label its vertices in cyclic order as $v_1,\dots,v_g$.  Let $f_i=\deg(v_i)-2$
be the number of leaves adjacent to $v_i$.  Note that $G$ can be specified up to
isomorphism by the cyclically ordered list of numbers $f_1,\dots,f_g$.

Let $L$ denote the set of leaf edges of $G$.  The subsets of $E(G)$ whose
type is singleton-free are precisely those that contain $L$ as a subset.
In particular, $L$ itself is the unique edge set whose type is 
a singleton-free partition of length~$g$.  Thus $\type(L)=\{f_1-1,f_2-1,
\dots,f_g-1\}$ can be recovered from $\CSF_G$.  Moreover, there are
precisely $g$ edge sets whose type is a singleton-free partition
of length~$g-1$; these edge sets are of the form $L\cup\{e\}$ for some
$e\in C$.  Just as in Theorem~\ref{distinct-caterpillar-thm},
the types of these edge sets are all distinct, and they specify which
pairs of the $f_i$ correspond to adjacent vertices of the cycle.
\end{proof}

%=================================================================
\bibliographystyle{amsalpha}

\end{document}